\documentclass{article}
\usepackage{latexsym}
\usepackage{amsfonts}
\newtheorem{theorem}{Theorem}

\newtheorem{proposition}{Proposition}
\newtheorem{claim}{Claim}

\begin{document}
\def\Var{{\rm Var}}
\def\grpz{{\mathbb Z}}
\def\grpc{{\mathbb C}}
\def\fbar{{\overline f}}
\def\qed{{\hfill $\Box$}}
\title{On the Chung-Diaconis-Graham random process}
\author{Martin Hildebrand\\Department of Mathematics and Statistics\\
University at Albany\\State University of New York\\Albany, NY 12222}

\maketitle

\begin{abstract}
Chung, Diaconis, and Graham considered random processes of the form
$X_{n+1}=2X_n+b_n \pmod p$ where $X_0=0$, $p$ is odd, and $b_n$ for
$n=0,1,2,\dots$ are i.i.d. random variables on $\{-1,0,1\}$.
If $\Pr(b_n=-1)=\Pr(b_n=1)=\beta$ and $\Pr(b_n=0)=1-2\beta$, they
asked which value of $\beta$ makes $X_n$ get close to uniformly distributed 
on the integers mod $p$ the slowest.  In this paper, we
extend the results of Chung, Diaconis, and Graham in the case
$p=2^t-1$ to show that for $0<\beta\le 1/2$, there is no such
value of $\beta$.
\end{abstract}

\section{Introduction}
In \cite{cdg}, Chung, Diaconis, and Graham considered random processes of
the form $X_{n+1}=2X_n+b_n\pmod p$ where $p$ is an odd integer,
$X_0=0$, and $b_0, b_1, b_2, \dots$ are i.i.d. random variables.
This process is also described in Diaconis~\cite{diaconis}, and
generalizations involving random processes of the form
$X_{n+1}=a_nX_n+b_n \pmod p$ where $(a_i,b_i)$ for $i=0,1,2,\dots$ are i.i.d.
were considered by the author in \cite{mvh} and \cite{mvhima}.
A question asked in \cite{cdg} 
concerns cases where $\Pr(b_n=1)=\Pr(b_n=-1)=\beta$
and $\Pr(b_n=0)=1-2\beta$. If $\beta=1/4$ or $\beta=1/2$, then $P_n$ is
close to the uniform distribution (in variation distance)
 on the integers mod $p$ if $n$ is a large enough multiple of
$\log p$ where $P_n(s)=\Pr(X_n=s)$. If $\beta=1/3$, however, for $n$
a small enough multiple of $(\log p)\log(\log p)$, the variation
distance $\|P_n-U\|$ is far from $0$ for certain values of $p$
such as $p=2^t-1$. Chung, Diaconis, and Graham comment ``It would be 
interesting to know which value of $\beta$ maximizes the value of $N$
required for $\|P_N-U\|\rightarrow 0$.''

If $\beta=0$, then $X_n=0$ with probability $1$ for all $n$. Thus we shall
only consider the case $\beta>0$. We shall show that unless $\beta=1/4$
or $\beta=1/2$, then there exists a value $c_{\beta}>0$ such that
for certain values of $p$ (namely $p=2^t-1$), if
$n\le c_{\beta}(\log p)\log(\log p)$, then $\|P_n-U\|\rightarrow 1$ 
as $t\rightarrow\infty$.
Furthermore, one can have $c_{\beta}\rightarrow \infty$ as $\beta
\rightarrow 0^{+}$. Work of the author~\cite{mvh} shows that for each $\beta$,
there is a value $c_{\beta}^{\prime}$ such that
if $n\ge c_{\beta}^{\prime}(\log p)\log(\log p)$, then $\|P_n-U\|
\rightarrow 0$ as $p\rightarrow\infty$. Thus one may conclude 
that there is no value of $\beta$
which maximizes the value of $N$ required for $\|P_N-U\|\rightarrow 0$.

This paper will consider a broader class of distributions for $b_n$. In
particular, $\Pr(b_n=1)$ need not equal $\Pr(b_n=-1)$.
The main argument here relies on a generalization of an 
argument in \cite{cdg}.

\section{Notation and Main Theorem}

Recall that the variation distance of a probability $P$ on 
a finite group $G$ from the uniform distribution on $G$ is given by
\begin{eqnarray*}
\|P-U\|&=&{1\over 2}\sum_{s\in G}\left|P(s)-1/|G|\right| \\
&=&\max_{A\subseteq G}\left|P(A)-U(A)\right| \\
&=&\sum_{s:P(s)>1/|G|}\left|P(s)-1/|G|\right|
\end{eqnarray*}

The following assumptions are used in the main theorem.
Suppose $\Pr(b_n=1)=a$, $\Pr(b_n=0)=b$, and $\Pr(b_n=-1)=c$. We assume
$a+b+c=1$ and $a$, $b$, and $c$ are all less than $1$. Suppose
$b_0, b_1, b_2, \dots$ are i.i.d. and $X_0=0$. Suppose
$X_{n+1}=2X_n+b_n \pmod p$ and $p$ is odd. Let $P_n(s)=\Pr(X_n=s)$.
The theorem itself follows:
\begin{theorem}
\label{mainthm}
Case 1: Suppose either $b=0$ and $a=c=1/2$ or $b=1/2$. If
$n>c_1\log_2 p$ where $c_1>1$ is constant, then
$\|P_n-U\|\rightarrow 0$ as $p\rightarrow\infty$ where $p$ is an
odd integer.

Case 2: Suppose $a$, $b$, and $c$ do not satisfy the conditions in
Case 1. Then there exists a value $c_2$ (depending on $a$, $b$, and 
$c$) such that if $n<c_2(\log p)\log(\log p)$ and $p=2^t-1$,
then $\|P_n-U\|\rightarrow 1$ as $t\rightarrow\infty$.
\end{theorem}

\section{Proof of Case 1}

First let's consider the case where $b=1/2$. Then $b_n=e_n+d_n$ where $e_n$
and $d_n$ are independent random variables with $\Pr(e_n=0)=\Pr(e_n=1)=1/2$,
$\Pr(d_n=-1)=2c$, and $\Pr(d_n=0)=2a$. (Note that here $a+c=1/2=b$. Thus 
$2a+2c=1$.) Observe that
\begin{eqnarray*}
X_n&=&\sum_{j=0}^{n-1}2^{n-1-j}b_j\pmod p \\
&=&\sum_{j=0}^{n-1}2^{n-1-j}e_j+\sum_{j=0}^{n-1}2^{n-1-j}d_j\pmod p
\end{eqnarray*}
Let
\[
Y_n=\sum_{j=0}^{n-1}2^{n-1-j}e_j\pmod p.
\]
If $P_n$ is the probability distribution of $X_n$ (i.e.
$P_n(s)=\Pr(X_n=s)$) and $Q_n$ is the probability distribution of $Y_n$, 
then the independence of $e_n$ and $d_n$ implies $\|P_n-U\|\le
\|Q_n-U\|$.
Observe that on the integers,
$\sum_{j=0}^{n-1}2^{n-1-j}e_j$ is uniformly distributed on the set
$\{0,1,\dots,2^n-1\}$. Each element of the integers mod $p$ appears
either $\lfloor 2^n/p\rfloor$ times or $\lceil 2^n/p\rceil$ times. Thus
\[
\|Q_n-U\|\le p\left({\lceil 2^n/p\rceil\over 2^n}-{1\over p}\right)\le {p\over 
2^n}.
\]
If $n>c_1\log_2p$ where $c_1>1$, then $2^n>p^{c_1}$ and $\|Q_n-U\|\le
1/p^{c_1-1}\rightarrow 0$ as $p\rightarrow\infty$.

The case where $b=0$ and $a=c=1/2$ is alluded to in \cite{cdg} and
left as an exercise.\qed

\section{Proof of Case 2}

The proof of this case follows the proof of Theorem 2 in \cite{cdg}
with some modifications. First we assume $b>0$.

Define, as in \cite{cdg}, the separating function
$f: \grpz/p\grpz\rightarrow\grpc$ by
\[
f(k):=\sum_{j=0}^{t-1}q^{k2^j}
\]
where $q:=q(p):=e^{2\pi i/p}$. We shall suppose $n=rt$ where
$r$ is an integer of the form $r=\delta\log t-d$ for a fixed value
$\delta$.

If $0\le j\le t-1$, define
\[
\Pi_j:=\prod_{\alpha=0}^{t-1}\left(aq^{(2^{\alpha}(2^j-1))}+b+
cq^{-(2^{\alpha}(2^j-1))}\right).
\]
Note that if $a=b=c=1/3$, then this expression is the same as $\Pi_j$
defined in the proof of Theorem 2 in \cite{cdg}.

As in the proof of Theorem 2 in \cite{cdg},
$E_U(f)=0$ and $E_U(f{\fbar})=t$. Furthermore
\begin{eqnarray*}
E_{P_n}(f)&=&\sum_k P_n(k)f(k)\\
&=&\sum_k \sum_{j=0}^{t-1}P_n(k)q^{k2^j}\\
&=&\sum_{j=0}^{t-1}\hat P_n(2^j)\\
&=&\sum_{j=0}^{t-1}\prod_{\alpha=0}^{t-1}\left(
aq^{2^{\alpha}2^j/p}+b+cq^{-2^{\alpha}2^j/p}\right)^r\\
&=&t\Pi_1^r.
\end{eqnarray*}

Also note
\begin{eqnarray*}
E_{P_n}(f{\fbar})&=&\sum_kP_n(k)f(k)\fbar(k)\\
&=&\sum_k\sum_{j,j^{\prime}}P_n(k)q^{k(2^j-2^{j^{\prime}})}\\
&=&\sum_{j,j^{\prime}}\hat P_n(2^j-2^{j^{\prime}})\\
&=&\sum_{j,j^{\prime}}\prod_{\alpha=0}^{t-1}\left(aq^{2^{\alpha}
(2^j-2^{j^{\prime}})}+b+cq^{-2^{\alpha}(2^j-2^{j^{\prime}})}\right)^r\\
&=&t\sum_{j=0}^{t-1}\Pi_j^r.
\end{eqnarray*}
(Note that the expressions for $E_{P_N}(f)$ and $E_{P_N}(f\fbar)$
in the proof of Theorem 2 of \cite{cdg} have some minor misprints.)

The (complex) variances of $f$ under $U$ and $P_n$ are 
$\Var_U(f)=t$ and
\begin{eqnarray*}
\Var_{P_n}(f)&=&E_{P_n}(|f-E_{P_n}(f)|^2)\\
&=&E_{P_N}(f\fbar)-E_{P_n}(f)E_{P_n}(\fbar)\\
&=&t\sum_{j=0}^{t-1}\Pi_j^r-t^2|\Pi_1|^{2r}.
\end{eqnarray*}
Like \cite{cdg}, we use the following complex form of Chebyshev's
inequality for any $Q$:
\[
Q\left(\left\{x:|f(x)-E_Q(f)|\ge\alpha\sqrt{\Var_Q(f)}\right\}\right)\le 
1/\alpha^2
\]
where $\alpha>0$.
Thus
\[
U\left(\left\{x:|f(x)|\ge\alpha t^{1/2}\right\}\right)\le 1/\alpha^2
\]
and
\[
P_n\left(\left\{x:|f(x)-t\Pi_1^r|\ge \beta
\left(t\sum_{j=0}^{t-1}\Pi_j^r-t^2|\Pi_1|^{2r}\right)^{1/2}\right\}\right)
\le 1/\beta^2.
\]
Let $A$ and $B$ denote the complements of these 2 sets; thus $U(A)\ge 
1-1/\alpha^2$ and $P_n(B)\ge 1-1/\beta^2$. If $A$ and $B$ are disjoint, 
then $\|P_n-U\|\ge1-1/\alpha^2-1/\beta^2$.

Suppose $r$ is an integer with
\[
r={\log t\over 2 \log(1/|\Pi_1|)}-\lambda
\]
where $\lambda\rightarrow\infty$ as $t\rightarrow\infty$ but $\lambda\ll
\log t$. Then $t|\Pi_1|^r=t^{1/2}|\Pi_1|^{-\lambda}\gg t^{1/2}$.
Observe that the fact $a$, $b$, and $c$ do not satisfy the conditions in
Case 1 implies $|\Pi_1|$ is bounded away from $0$ as $t\rightarrow\infty$.
Furthermore $|\Pi_1|$ is bounded away from $1$ for a given $a$, $b$, and
$c$.

In contrast, let's consider what happens to $|\Pi_1|$ if $a$, $b$, and $c$
do satisfy the condition in Case 1. If $b=1/2$, then the $\alpha=t-1$ term
in the definition of $\Pi_1$ converges to $0$ as $t\rightarrow\infty$
and thus $\Pi_1$ also converges to $0$ as $t\rightarrow\infty$ since each other
term has length at most $1$. If $a=c=1/2$ and $b=0$, then the
$\alpha=t-2$ term in the definition of $\Pi_1$ converges to $0$
as $t\rightarrow\infty$ and thus $\Pi_1$ also converges to $0$
as $t\rightarrow\infty$.

\begin{claim}
\[
{1\over t}\sum_{j=0}^{t-1}\left({\Pi_j\over |\Pi_1|^2}\right)^r\rightarrow 1
\]
as $t\rightarrow\infty$.
\end{claim}

Note that this claim implies $(\Var_{P_n}(f))^{1/2}=o(E_{P_n}(f))$ and thus
Case 2 of Theorem~\ref{mainthm} follows.

Note that $\Pi_0=1$. By Proposition~\ref{antisym} below, ${\overline \Pi_j}=
\Pi_{t-j}$. Thus $t\sum_{j=0}^{t-1}\Pi_j^r$ is real. Also note that
since $\Var_{P_n}(f)\ge 0$, we have
\[
{t\sum_{j=0}^{t-1}\Pi_j^r \over t^2|\Pi_1|^{2r}}\ge 1.
\]
Thus to prove the claim, it suffices to show
\[
{1\over t}\sum_{j=0}^{t-1}\left(|\Pi_j|\over|\Pi_1|^2\right)^r\rightarrow 1.
\]

\begin{proposition}
\label{antisym}
${\overline \Pi_j}=\Pi_{t-j}$.
\end{proposition}

{\it Proof:} Note that
\[
{\overline \Pi_j}=\prod_{\alpha=0}^{t-1}\left(aq^{-(2^{\alpha}(2^j-1))}
+b+cq^{(2^{\alpha}(2^j-1))}\right)
\]
and
\[
\Pi_{t-j}=\prod_{\beta=0}^{t-1}\left(aq^{(2^{\beta}(2^{t-j}-1))}+b+
cq^{-(2^{\beta}(2^{t-j}-1))}\right).
\]

If $j\le \beta\le t-1$, then note
\begin{eqnarray*}
2^{\beta}(2^{t-j}-1)&=&2^{\beta-j}(2^t-2^j)\\
&=&2^{\beta-j}(1-2^j)\pmod p\\
&=&-2^{\beta-j}(2^j-1).
\end{eqnarray*}
Thus the terms in $\Pi_{t-j}$ with $j\le\beta\le t-1$
are equal to the terms in ${\overline \Pi_j}$ with $0\le\alpha\le t-j-1$.

If $0\le\beta\le j-1$, then note
\begin{eqnarray*}
2^{\beta}(2^{t-j}-1)&=&2^{t+\beta}(2^{t-j}-1)\pmod p\\
&=&2^{t+\beta-j}(2^t-2^j)\\
&=&2^{t+\beta-j}(1-2^j)\pmod p\\
&=&-2^{t+\beta-j}(2^j-1).
\end{eqnarray*}
Thus the terms in $\Pi_{t-j}$ with $0\le\beta\le j-1$ are equal to the
terms in ${\overline\Pi_j}$ with $t-j\le\alpha\le t-1$.\qed

Now let's prove the claim. Let $G(x)=|ae^{2\pi ix}+b+ce^{-2\pi ix}|$. Thus
\[|\Pi_j|=\prod_{\alpha=0}^{t-1}G(2^{\alpha}(2^j-1)/p).\]
Note that if $0\le x<y\le 1/4$, then $G(x)>G(y)$.
On the interval $[1/4,1/2]$, where $G$ increases and where $G$ decreases
depends on $a$, $b$, and $c$.

We shall prove a couple of facts analogous to facts in \cite{cdg}.

{\it Fact 1:} There exists a value $t_0$ (possibly depending on
$a$, $b$, and $c$) such that if $t>t_0$, then 
$|\Pi_j|\le|\Pi_1|$ for all $j\ge 1$.

Since $G(x)=G(1-x)$, in proving this fact we may assume without loss
of generality that $2\le j\le t/2$.
Note that
\[
|\Pi_j|=\prod_{i=0}^{t-j-1}G\left({2^{i+j}-2^i\over p}\right)
\prod_{i=0}^{j-1}G\left({2^{i+t-j}-2^i\over p}\right).
\]
We associate factors $x$ from $|\Pi_j|$ with corresponding factors $\pi(x)$
of $|\Pi_1|$ in a manner similar to that in \cite{cdg}. For $0\le i\le t-j-2$,
associate $G((2^{i+j}-2^i)/p)$ with $G(2^{i+j-1}/p)$. Note that for 
$0\le i\le t-j-2$, we have $G((2^{i+j}-2^i)/p)\le G(2^{i+j-1}/p)$.
For $0\le i\le j-3$, associate $G((2^{i+t-j}-2^i)/p)$ in $|\Pi_j|$
with $G(2^i/p)$ in $|\Pi_1|$. Note that for $0\le i\le j-3$, we have
$G((2^{i+t-j}-2^i)/p)\le G(2^i/p)$.

The remaining terms in $|\Pi_j|$ are
\[
G\left({2^{t-1}-2^{t-j-1}\over p}\right)G\left({2^{t-1}-2^{j-1}\over p}\right)
G\left({2^{t-2}-2^{j-2}\over p}\right)
\]
and the remaining terms in $|\Pi_1|$ are
\[
G\left({2^{t-1}\over p}\right)G\left({2^{t-2}\over p}\right)G\left(
{2^{j-2}\over p}\right).
\]
It can be shown that since $b>0$,
\[
\limsup_{t\rightarrow\infty}{
G\left({2^{t-1}-2^{t-j-1}\over p}\right)G\left({2^{t-1}-2^{j-1}\over p}\right)
G\left({2^{t-2}-2^{j-2}\over p}\right)
\over
G\left({2^{t-1}\over p}\right)G\left({2^{t-2}\over p}\right)G\left(
{2^{j-2}\over p}\right)
}\le{M\over G(0)}<1\]
where $M=\sup_{x\in[1/4,1/2]}G(x)$. Note $G(0)=1$. To see the fact that
$M<1$ if $b>0$, consider $H(x)=|ae^{2\pi ix}+b|$ if $a>0$. $H(x)$ is decreasing
on $[0,1/2]$, and $a+b-H(1/4)>0$. For $x\in[1/4,1/2]$, we have
$G(x)\le H(x)+c=1-(a+b-H(x))\le 1-(a+b-H(1/4))$. If $a=0$, then $c>0$
and a similar argument involving $|ce^{-2\pi ix}+b|$ applies.
Indeed, for some $t_0$, if $t>t_0$ and $2\le j\le t/2$,
\begin{eqnarray*}
&&G\left({2^{t-1}-2^{t-j-1}\over p}\right)G\left({2^{t-1}-2^{j-1}\over p}\right)
G\left({2^{t-2}-2^{j-2}\over p}\right)
\\&\le&
G\left({2^{t-1}\over p}\right)G\left({2^{t-2}\over p}\right)G\left(
{2^{j-2}\over p}\right).
\end{eqnarray*}\qed

{\it Fact 2:} There exists a value $t_1$ (possibly depending on $a$, $b$,
and $c$) such that if $t>t_1$, then the following holds. There is a constant
$c_0$ such that for $t^{1/3}\le j\le t/2$, we have
\[
{|\Pi_j|\over |\Pi_1|^2}\le1+{c_0\over 2^j}
\]

To prove this fact, we associate, for $i=0,1,.\dots,j-1$, the terms
\[G\left({2^{t-i-1}-2^{j-i-1}\over p}\right)
G\left({2^{t-i-1}-2^{t-j-i-1}\over p}\right)\] in $|\Pi_j|$
with the terms \[\left(G\left({2^{t-i-1}\over p}\right)
\right)^2\] in $|\Pi_1|^2$. Suppose
$A=\max|G^{\prime}(x)|$. Note that $A<\infty$. Then
\[
\left|G\left({2^{t-i-1}-2^{j-i-1}\over p}\right)\right|
\le
\left|G\left({2^{t-i-1}\over p}\right)\right|+A{2^{j-i-1}\over p}.
\]
Thus
\[
{
\left|G\left({2^{t-i-1}-2^{j-i-1}\over p}\right)\right|
\over
\left|G\left({2^{t-i-1}\over p}\right)\right|
}\le 1+A{2^{j-i-1}\over p\left|G\left({2^{t-i-1}\over p}\right)\right|}.
\]
Likewise
\[
{
\left|G\left({2^{t-i-1}-2^{t-j-i-1}\over p}\right)\right|
\over
\left|G\left({2^{t-i-1}\over p}\right)\right|
}\le 1+A{2^{t-j-i-1}\over p\left|G\left({2^{t-i-1}\over p}\right)\right|}.
\]
Since we do not have the conditions for Case 1, there is a positive value
$B$ and value $t_2$ such that if $t>t_2$, then
$|G(2^{t-i-1}/p)|>B$ for all $i$ with $0\le i\le j-1$. By an exercise, one
can verify 
\[
\prod_{i=0}^{j-1}
{
\left|G\left({2^{t-i-1}-2^{j-i-1}\over p}\right)G\left({2^{t-i-1}-2^{t-j-i-1}
\over p}\right)\right|
\over
\left|G\left({2^{t-i-1}\over p}\right)\right|^2
}\le
1+{c_3\over 2^j}
\]
for some value $c_3$ not depending on $j$.

Note that the remaining terms in $|\Pi_j|$ all have length less than $1$. The
remaining terms in $|\Pi_1|^2$ are
\[
\prod_{i=j}^{t-1}\left|G\left({2^{t-i-1}\over p}\right)\right|^2.
\]
Since $G^{\prime}(0)=0$, there are positive constants $c_4$ and $c_5$ such 
that \[\left|G\left({2^{t-i-1}\over p}\right)\right|\ge 
1-c_4\left({2^{t-i-1}\over p}\right)^2\ge\exp\left(-c_5{2^{t-i-1}\over p}\right)
\]
for $i\ge j\ge t^{1/3}$.
Observe
\begin{eqnarray*}
\prod_{i=j}^{t-1}\exp\left(-c_5{2^{t-i-1}\over p}\right)&=&
\exp\left(-c_5\sum_{i=j}^{t-1}2^{t-i-1}/p\right)\\
&=&\exp\left(-c_5\sum_{k=0}^{t-j-1}2^k/p\right)\\
&=&\exp\left(-c_5{2^{t-j}-1\over 2^t-1}\right)\\
&>&\exp\left(-c_5{2^{t-j}\over 2^t}\right)\\
&=&\exp(-c_5/2^j)>1-c_5/2^j.
\end{eqnarray*}

There exists a constant $c_0$ such that
\[
{1+c_3/2^j \over (1-c_5/2^j)^2}\le 1+c_0/2^j
\]
for $j\ge 1$.

Thus, as in \cite{cdg},
\[
\sum_{t^{1/3}\le j\le t/2}\left|\left({|\Pi_j|\over |\Pi_1|^2}\right)^r-1
\right|\le {c_6tr\over 2^{t^{1/3}}}<{c_7\over 2^{t^{1/4}}}
\]
for values $c_6$ and $c_7$. Since $|\Pi_j|=|\Pi_{t-j}|$,
\begin{eqnarray*}
{1\over t}\sum_{j=0}^{t-1}\left({|\Pi_j|\over |\Pi_1|^2}\right)^r&\le&
{1\over t}{1 \over |\Pi_1|^{2r}}+{2\over t}\left(\sum_{1\le j< t^{1/3}}
\left({|\Pi_j|\over |\Pi_1|^2}\right)^r+\sum_{t^{1/3}\le j\le t/2}
\left({|\Pi_j|\over |\Pi_1|^2}\right)^r\right)\\
&=&
1+o(1)
\end{eqnarray*}
as $t\rightarrow\infty$. Thus Fact 2, the claim, and Theorem~\ref{mainthm} if 
$b>0$ are
proved.

If $b=0$, Theorem~\ref{mainthm} can be proved by considering the following
random processes on the integers mod $p$:

1. $X_0=0$ and $X_{n+1}=2X_n+b_n \pmod p$ where $P(b_n=1)=a$ and $P(b_n=-1)=1-a$

2. $Y_0=0$ and $Y_{n+1}=2Y_n+d_n \pmod p$ where $P(d_n=2)=a$ and $P(d_n=0)=1-a$

3. $Z_0=0$ and $Z_{n+1}=2Z_n+e_n \pmod p$ where $P(e_n=1)=a$ and $P(e_n=0)=1-a$

If $P_n(s)=Pr(X_n=s)$, $Q_n(s)=Pr(Y_n=s)$, and $R_n(s)=Pr(Z_n=s)$, then
$\|P_n-U\|=\|Q_n-U\|=\|R_n-U\|$. To see this, we can let $d_n=2e_n$ and
$b_n=d_n-1$ so that $Y_n=2Z_n$ and
$X_n=Y_n-\sum_{j=0}^{n-1}2^j$ for $n\ge 1$. To conclude the case where
$b=0$, use the argument with $b\ne 0$ on the third random process (provided
that $a\ne 1/2$ so that we are in Case 2).

\qed

The next proposition considers what happens as we vary the values
$a$, $b$, and $c$.

\begin{proposition}
\label{smallbeta}
If $a=c=\beta$ and $b=1-2\beta$ and $m_{\beta}=\liminf_{t\rightarrow
\infty}|\Pi_1|$, then $\lim_{\beta\rightarrow 0^{+}}m_{\beta}=1$.
\end{proposition}

{\it Proof:} Suppose $\beta<1/4$. Then
\[
\Pi_1=\prod_{\alpha=0}^{t-1}\left((1-2\beta)+2\beta\cos(2\pi2^{\alpha}/p)
\right).
\]
Let $h(\alpha)=(1-2\beta)+2\beta\cos(2\pi2^{\alpha}/p)$. Note that
\begin{eqnarray*}
\lim_{\beta\rightarrow 0^{+}}h(t-1)&=&1\\
\lim_{\beta\rightarrow 0^{+}}h(t-2)&=&1\\
\lim_{\beta\rightarrow 0^{+}}h(t-3)&=&1
\end{eqnarray*}
Furthermore, for some constant $\gamma>0$, one can show
\[
h(\alpha)>\exp(-\beta\gamma(2^{\alpha}/p)^2)
\]
if $2^{\alpha}/p\le 1/8$ and $0<\beta<1/10$.
So
\begin{eqnarray*}
\prod_{\alpha=0}^{t-4}h(\alpha)&>&\prod_{\alpha=0}^{t-4}\exp(-\beta
\gamma(2^{\alpha}/p)^2)
\\
&=&\exp\left(-\beta\gamma\sum_{\alpha=0}^{t-4}(2^{\alpha}/p)^2\right)\\
&>&\exp(-\beta\gamma2^{2(t-4)}(4/3)/p^2)\rightarrow 1
\end{eqnarray*}
as $\beta\rightarrow 0^{+}$.\qed

Recalling that
\[
r={\log t\over 2\log(1/|\Pi_1|)}-\lambda,
\]
we see that $1/(2\log(1/|\Pi_1|))$ can be made arbitrarily large by choosing
$\beta$ small enough. Thus there exist values $c_{\beta}\rightarrow\infty$ as
$\beta\rightarrow 0^{+}$ such that if $n\le c_{\beta}(\log p)\log(\log p)$,
then
$\|P_n-U\|\rightarrow 1$ as $t\rightarrow\infty$.

\section{Problems for further study}

One possible problem is to see if in some sense, there is a value of
$\beta$ on $[1/4,1/2]$ which maximizes the value of $N$ required for
$\|P_N-U\|\rightarrow 0$; to consider such a question, one might
restrict $p$ to values such that $p=2^t-1$.

Another possible question considers the behavior of these random processes
for almost all odd $p$. For $\beta=1/3$, Chung, Diaconis, and Graham~\cite{cdg}
showed that a multiple of $\log p$ steps suffice for almost all odd $p$.
While their arguments should be adaptable with the change of appropriate 
constants to a broad range of choices of $a$, $b$, and $c$ in Case 2,
a more challenging question is to determine for which $a$, $b$, and $c$
in Case 2 (if any), $(1+o(1))\log_2p$ steps suffice for almost all odd $p$.

\section{Acknowledgments}

The author thanks Ron Graham for mentioning the problem at a January, 2005,
conference on the Mathematics of Persi Diaconis. The author also thanks Robin 
Pemantle for a conversation on the topic and the participants in the
Probability and Related Fields Seminar at the University at Albany for
listening to some ideas on the problem.
The author also thanks Ravi Montenegro for pointing out a flaw in the original
argument, which appeared in \cite{mvhecp}.


\begin{thebibliography}{99}

\bibitem{cdg} Chung, F., Diaconis, P., and Graham, R. A random walk
problem arising in random number generation. {\it Ann. Probab.} {\bf 15}
(1987), 1148-1165.

\bibitem{diaconis} Diaconis, P. {\it Group Representations in
Probability and Statistics}. Institute of Mathematical Statistics, 1988.

\bibitem{mvh} Hildebrand, M. Random processes of the form
$X_{n_+1}=a_nX_n+b_n\pmod p$. {\it Ann. Probab.} {\bf 21} (1993), 710-720.

\bibitem{mvhima} Hildebrand, M. Random processes of the form
$X_{n+1}=a_nX_n+b_n \pmod p$ where $b_n$ takes on a single value.
pp. 153-174, {\it Random Discrete Structures}, ed. Aldous and Pemantle.
Springer-Verlag, 1996.

\bibitem{mvhecp} Hildebrand, M. On the Chung-Diaconis-Graham random process.
{\it Electron. Comm. in Probab.} {\bf 11} (2006), 347-356.
\end{thebibliography}
\end{document}